\documentclass[12pt]{article}
\usepackage{amssymb,latexsym}
\usepackage{amsmath,bbm,amsthm}
\usepackage{latexsym}
\usepackage{euscript}
\usepackage[active]{srcltx}

\oddsidemargin=-0.1cm \tolerance 9000 \theoremstyle{definition}

\def\finproof{\hfill\hbox{\vrule width1.0ex height1.5ex}\vspace{2mm}}

\begin{document}

\begin{center}

{{
\textbf{ \sc
Local properties of solutions to non-autonomous parabolic PDEs 
with state-dependent delays}}} 

\vskip7mm

{
\textsc{Alexander V. Rezounenko }}

\smallskip

Department of Mechanics and Mathematics \\ Kharkiv National
University, 4,
Svobody Sqr., Kharkiv, 61077, Ukraine \\ 
rezounenko@univer.kharkov.ua

\end{center}

\bigskip

{\bf Abstract.} A wide class of non-autonomous nonlinear parabolic
partial differential equations with delay is studied. We allow in
our investigations different types of delays such as constant,
time-dependent, state-dependent (both discrete and distributed) to
be presented simultaneously. The main difficulties arise due to the
presence of discrete state-dependent delays since the nonlinear
delay term is not Lipschitz on the space of continuous functions. We
find conditions for the local existence, uniqueness and study the
invariance principle.




\section{Introduction}
\marginpar{\tiny March 22, 2011}

We consider non-autonomous parabolic partial differential equations
(PDEs) with delay. Studying of this type of equations is based on
1) the well-developed theories of the delayed ordinary differential equations
(ODEs) 
\cite{Hale,Hale_book,Walther_book} 
 and 2) PDEs without delays \cite{Hadamard-1902,Hadamard-1932,Lions-Magenes-book,Lions}. Under
certain assumptions both types of equations describe a kind of
dynamical systems that are infinite-dimensional, see
\cite{Babin-Vishik,Temam_book,Chueshov_book} and references therein;
see also
\cite{travis_webb,Fitzgibbon-JDE-1978,Martin-Smith-TAMS-1990,Chueshov-JSM-1992,Cras-1995}
and the monograph \cite{Wu_book} that are close to our work.

In 
evolution systems arising in applications the presented delays are
frequently state-dependent (SDDs). The theory of such equations,
especially the ODEs, is rapidly developing and many deep results
have been obtained up to now (see e.g.
\cite{Mallet-Paret,Walther_JDE-2003,Walther_JDDE-2007,Krisztin-2003}
and also the survey paper \cite{Hartung-Krisztin-Walther-Wu-2006}
for details and references). The PDEs with state-dependent delays
were first studied in
\cite{Rezounenko-Wu-2006,Hernandez-2006,Rezounenko_JMAA-2007}. An
alternative approach to the PDEs with discrete SDDs is proposed in
\cite{Rezounenko_NA-2009}. Approaches to equations with discrete and
distributed SDDs are different. Even in the case of ODEs, the
discrete SDD essentially complicates the study since, in general,
the corresponding {\it nonlinearity} is {\it not} locally Lipschitz continuous
on open subsets of the space of continuous functions, and familiar
results on existence, uniqueness, and dependence of solutions on
initial data and parameters from, say \cite{Hale_book,Walther_book}
fail (see \cite{Winston-1970} for an example of the non-uniqueness
and  \cite{Hartung-Krisztin-Walther-Wu-2006} for more details). It is important to mention
that due to the discrete SDDs such equations are inherently {\it nonlinear}.
 In
this work, in contrast to previous investigations, we consider a
model where two different types of SDDs (discrete and distributed)
are presented simultaneously (by Stieltjes integral).
Moreover, all the
assumptions on the delay (see (A1)-(A5) below) allow the dynamics
when along a solution the number and values of discrete SDDs may
change, the whole discrete and/or distributed delays may vanish,
disappear and appear again. This property makes it possible 
to study "flexible" models
where some subsets of the phase space are described by equations
with purely discrete SDDs, and others by equations with purely
distributed SDDs, and there are subsets which need the general
(combined) type of the delay. A solution could be in different
subsets at different time moments. This property particularly means
that not only the values of the delays are state-dependent, but the
{\it type} of the delay is {\it state-dependent} as well.

The first goal of the present paper is to study the basic properties
of solutions - the existence and uniqueness as well as to extend the
fundamental invariance principle to the case of PDEs with discrete
SDDs. The second goal is to attract attention of researchers from
such fields as, for example, mathematical biology and physics to
this wide class of delay equations and emphasize that the crucial
assumption on the delay (see (A5) below) is an "inner property" of
the delay which could be successfully used in a wide range of other
delay systems. We hope that our results provide a basement for further study
of qualitative (asymptotic) properties of solutions.

The existence and uniqueness results for a particular case of
autonomous systems were announced in \cite{Rezounenko-CR-2011}. For
a survey of the existing literature on the invariance principle see
\cite{Ruess-TAMS-2009}. In the present paper the emphasis is on the
delayed term, not on the partial differential operator. To the best
of our knowledge the invariance principle for PDEs with SDDs has not
been studied before.

\section{Formulation of the model and examples}

Let $X$ be a Banach space with the norm $||\cdot ||$, let $r>0$ be a
constant. Denote by $C\equiv C([-r,0];X)$ the space of continuous
functions $\varphi : [-r,0] \to X$ with the supremum norm $||\cdot
||_C.$  As usually for delay equations \cite{Hale,Hale_book}, for
any real $a\le b, t\in [a,b]$ and any continuous function $u :
[a-r,b] \to X$, we denote by $u_t$ the element of $C$ defined by the
formula $u_t = u_t(\theta) \equiv u(t+\theta)$ for $\theta\in
[-r,0].$ Consider an infinitesimal generator $A$   of a (compact)
$C_0$ semigroup $\{ e^{-At}\}_{t\ge 0}\equiv \{ T(t)\}_{t\ge 0}$ on
$X$ 
satisfying $||T(t)||\le  e^{\omega t}$ for all $t\ge 0$, where
$\omega\in \mathbb{R}$ is a fixed constant.

 We are interested in the following non-autonomous parabolic partial differential equation with
state-dependent delays (SDD)
\begin{equation}\label{sdd13-1}
{d u(t)\over dt} +Au(t) = B\big(t,u_t \big),\quad t\ge a
\end{equation}
with the initial condition
\begin{equation}\label{sdd13-ic}
  u_a=u|_{[a-r,a]}=\varphi \in C\equiv C([-r,0];X). 
\end{equation}

The delay term $B : \mathbb{R}\times C \to X$ has the form
\begin{equation}\label{sdd13-2}
B\big(t,\psi \big)\equiv G\big( t,\psi(0), F(t,\psi)\big),
\end{equation}
where $G : \mathbb{R}\times X \times X \to X$ is a continuous
 mapping and the delay functional $F:
\mathbb{R}\times C \to X$ is presented by a Stieltjes integral
(simultaneously includes discrete and distributed SDDs)
\begin{equation}\label{sdd13-3}
 F(t,\psi) \equiv \int^0_{-r}  p\left( t,\psi(\theta) \right)
 \cdot d g(\theta, t,\psi), \quad p: \mathbb{R}\times X\to X.
\end{equation}

Assumptions on $g$ are formulated below (see (A1)-(A5)).

The class of equations described by (\ref{sdd13-1}),(\ref{sdd13-2}),(\ref{sdd13-3})
is very wide and includes many equations which were intensively studied during past decades.
Below we mention just two examples and refer the reader to \cite{Wu_book,Rezounenko-Wu-2006}
for more references and discussion.

\medskip

{\it Example 1}. 
 Let $X=L^2(\Omega)$, where
$\Omega \subset \mathbb{R}^{N}$ is  a smooth bounded domain.
Operator $A$ is a densely-defined self-adjoint positive linear
operator with domain $D(A)\subset L^2(\Omega )$ and compact
resolvent, which means that $A: D(A)\to L^2(\Omega )$ generates an
analytic semigroup. If we choose $G(t,u,v)=v-d u$ ($d\ge 0$ is a
constant), then $B\big(t,u_t \big)=F(t,u_t) -d u$ and equation
(\ref{sdd13-1}) reads as
\begin{equation}\label{sdd13-4} \frac{\partial }{\partial
t}u(t,x)+Au(t,x)+du(t,x) = \big( F(u_t) \big)(x), \end{equation}
with, for example,
\begin{equation}\label{sdd13-5}
\big( F(\psi) \big)(x)\equiv \int^0_{-r} \left\{ \int_\Omega p\left(
\psi (\theta,y) \right) f(x-y) dy \right\} \cdot d g(\theta, \psi)
,\quad x\in \Omega,
\end{equation}
\noindent where  $f: \Omega -\Omega \to \mathbb{R}$   is a bounded
measurable function, $p:\mathbb{R}\to \mathbb{R}$. This non-local
{\it autonomous} equation is studied in \cite{Rezounenko-CR-2011}.
It is clear that the
integral delay term given by (\ref{sdd13-5}) includes the cases: \\
a) purely discrete SDDs: $\big( F(\psi) \big)(x)= \sum_k \int_\Omega
p\left( \psi(-\eta_k(\psi),y) \right) f(x-y)\, dy$; \\
b) purely distributed SDD: $\big( F(\psi) \big)(x) = \int^0_{-r}
\left\{ \int_\Omega p\left( \psi(\theta,y) \right) f(x-y) dy
\right\} \cdot \xi(\theta, \psi)\, d\theta.$
These cases have been studied in
\cite{Rezounenko-Wu-2006,Rezounenko_JMAA-2007,Rezounenko_NA-2009}.

Similarly, one may consider {\it local} delay terms (discrete and/or
distributed SDD)
\begin{equation}\label{sdd13-6}
\big( F(\psi) \big)(x)\equiv \int^0_{-r}  p\left( \psi(\theta,x)
\right) \cdot d g(\theta, \psi) ,\quad x\in \Omega.
\end{equation}
The above type of equations includes the diffusive Nicholson's blowflies
equation (see e.g. \cite{So-Yang}) with state-dependent delays, i.e.
equation (\ref{sdd13-4}) where $-A$ is the Laplace operator with 
Dirichlet or Neumann boundary conditions, $\Omega\subset
\mathbb{R}^{N}$ is a bounded domain with a smooth boundary, the
nonlinear (birth) function $p$ is given by $p(w)=p_1\cdot we^{-w}$,
$p_1\in \mathbb{R}$.

\medskip

{\it Example 2} ({\tt reaction-diffusion system with delay}).
Suppose $\Omega\subset {\mathbb{R}^N}$ is a bounded region with a
smooth boundary $\partial\Omega$, $\partial_n$ is the outward normal
derivative on $\partial\Omega$, $\triangle$ is the Laplacian
operator on $\Omega$. Consider the system
\begin{equation}\label{sdd13-ex-2-a} \left\{
\begin{array}{ll}
 \frac{\partial u^i}{\partial t}(x,t)= d_i \triangle u^i(x,t) +
G^i(t,F^i(t,u^1_t(x,\cdot),\ldots,u^m_t(x,\cdot))), & t>a,\, x\in
\Omega, \\ \\
  \alpha^i(x)u^i(x,t)+\partial_n u^i(x,t)=0,& t>a,\, x\in \partial\Omega,
  \\ \\
  u^i(x,a+\theta)=\varphi^i(x,\theta), & \theta\in [-r,0],\, x\in \Omega, \\
\end{array}
\right.
\end{equation}
where $i=1,\ldots,m$. In (\ref{sdd13-ex-2-a}), $d_i \ge 0$ and
$d_i=0$ we agree that no boundary condition applies to $u^i$,
$\alpha^i\in C^{1+\alpha}(\partial\Omega),\alpha\in (0,1)$.
Functions $G^i: \mathbb{R}^2\to \mathbb{R}$ are locally Lipschitz
and the delay functionals $F^i: \mathbb{R}\times
C([-r,0];\mathbb{R}^m) \to \mathbb{R}$ are presented by Stieltjes
integral (simultaneously includes discrete and distributed SDDs)
similar to (\ref{sdd13-3}). The system (\ref{sdd13-ex-2-a}) could be
presented in the form (\ref{sdd13-1})-(\ref{sdd13-2}) as follows. We
set $X\equiv C(\bar{\Omega};\mathbb{R}^m)$ and (see e.g.
\cite[p.5]{Martin-Smith-JRAM-1991}, 
\cite{Martin-Smith-TAMS-1990} and references therein) let $A^0_i$ be
the operator defined by $A^0_iy_i=d_i\Delta y_i $ (or $A^0_iy_i=0$
if $d_i=0$) on the domain $D(A^0_i)\equiv \{ y_i\in
C^1(\bar{\Omega})\cap C^2(\Omega); \, \alpha^i y^i+\partial_n y^i=0
\hbox{ on }
\partial\Omega\}$ (or $D(A^0_i)\equiv C(\bar{\Omega}) $ if
$d_i=0$)). The operator $A_i$ is the closure of $A^0_i$ on
$C(\bar{\Omega})$ and $A\equiv (A_i)^m_{i=1}.$ We denote
$\{T_i(t)\}_{t\ge 0}$ the $C_0$-semigroup on $C(\bar{\Omega})$
generated by $A_i$ and $T\equiv (T_i)^m_{i=1}.$ It is
well-known~\cite{Martin-Smith-TAMS-1990} that $T$ is a
$C_0$-semigroup on $X=C(\bar{\Omega};\mathbb{R}^m)$ that is analytic
(and compact if all $d_i > 0$) and $A$ is its generator. The
system~(\ref{sdd13-ex-2-a}) was studied (without state-dependent
delays), for example, in \cite{Martin-Smith-TAMS-1990}.

As an application one can consider the $n$-species Lotka-Volterra
model of competition with diffusion and delays given by
\begin{equation}\label{sdd13-ex-2-b} \left\{
\begin{array}{ll}
 \frac{\partial u^i}{\partial t}(x,t)= d_i \triangle u^i(x,t) +
b_i u^i(t,x)\left[ 1- \sum^n_{j=1} c_{ij} \int^0_{-r}u^j
(x,t+\theta))\, d g_{ij}(\theta,u_t)\right], & 
x\in
\Omega, \\ \\
\partial_n u^i(x,t)=0,& 
x\in \partial\Omega,
  \\ \\
  u^i(x,a+\theta)=\varphi^i(x,\theta), \quad \theta\in [-r,0], & 
   \\
\end{array}
\right.
\end{equation}
where $b_i, c_{ij}$ are positive constants and $g_{ij}$ are
nondecreasing with respect to the first coordinate and
$g_{ij}(0,\cdot)-g_{ij}(-r,\cdot)=1$. Many interesting properties of
this system (autonomous and without state-dependent delays) were
discussed in \cite{Martin-Smith-JRAM-1991} (see also references
therein).

\medskip

The approach developed in the present article is applicable to more general classes of
equations of the form (\ref{sdd13-1}) with the nonlinearity $B$, for example,
as follows (c.f. (\ref{sdd13-2})) %
$$
B\big(t,u_t \big)\equiv G\big( t,F^1(t,u_t),\ldots, F^k(t,u_t)\big),
$$
with $F^i$ be as in (\ref{sdd13-3}). We formulate our results for
$B$ given by (\ref{sdd13-2}) for the simplicity of presentation and
motivated by  (\ref{sdd13-ex-2-b}).



\section{Local existence and uniqueness}

The following assumptions on the {\it time-} and {\it
state}-dependent delay  are generalizations to the non-autonomous
case of the ones proposed in \cite{Rezounenko-CR-2011}.

\smallskip

{\bf (A1)} {\it For any $(t,\varphi)\in \mathbb{R}\times C $, the
function $[-r,0] \owns g(\cdot, t,\varphi)  \to \mathbb{R}$
is of bounded variation on $[-r,0].$ The variation $V^0_{-r}g$ of
$g$ is {\tt uniformly bounded} i.e.

$\exists M_{Vg}>0\, :\, \forall (t,\varphi)\in \mathbb{R}\times C
\quad \Rightarrow \quad V^0_{-r}g (\cdot,t,\varphi) \le M_{Vg}.$
}%

\smallskip

It is well-known that any Lebesgue-Stieltjes measure (associated with $g$) 
may be split into a sum of three measures: discrete, absolutely
continuous and singular ones. We will denote the corresponding
splitting of $g$ as follows

\begin{equation}\label{sdd13-4a}
 g(\theta, t,\varphi) = g_d(\theta, t,\varphi) + g_{ac}(\theta, t,\varphi) +g_s(\theta,
 t,\varphi)=g_d(\theta, t,\varphi) + g_{c}(\theta, t,\varphi),
\end{equation}
where $g_d(\theta, t,\varphi)$ is a step-function, $g_{ac}(\theta,
t,\varphi)$ is absolutely continuous and $g_s(\theta, t,\varphi)$ is
singular continuous as functions of their first coordinates (see
\cite{Kolmogorov-Fomin} for more details) and we denote the
continuous part by $g_c\equiv g_{ac}+g_s$.

Our next assumptions are 

{\bf (A2)} {\it For any $\theta\in [-r,0],$ the function $g_{c}$
is continuous with respect to 
its second and third coordinates i.e.
$\forall \theta \in [-r,0], \quad \forall
(t,\varphi),(t^n,\varphi^n)\in \mathbb{R}\times C  : (t^n,\varphi^n)
\to (t,\varphi) \hbox{ in } \mathbb{R}\times C \, (n\to +\infty)
\Rightarrow
g_{c}(\theta,t^n,\varphi^n) \to
g_{c}(\theta,t,\varphi).$
}%

\medskip

{\bf (A3)} 
 {\it The step-function $g_d(\theta,t,\varphi)$ is continuous with
 respect to $(t,\varphi)$ 
 in the sense that
 discontinuities of $g_d(\theta,t,\varphi)$ at points $\{ \theta_k\}\subset
 [-r,0]$ satisfy the property:
  there are {\tt continuous} functions  $\eta_k : \mathbb{R}\times C\to [0,r]$ and
  $h_k : \mathbb{R}\times C\to \mathbb{R}$ such that $\theta_k=-\eta_k(t,\varphi)$
 and $h_k (t,\varphi)$ is the jump of $g_d$ at point
$\theta_k=-\eta_k(t,\varphi)$ i.e $h_k(t,\varphi)\equiv g_d(
\theta_k+0,t,\varphi) - g_d( \theta_k-0,t,\varphi)$.

Taking into account that $g_d$ may, in general, have infinite
(countable) number of points of discontinuity $\{ \theta_k\}$, we
assume that the series $\sum_k  h_k (t,\varphi) $ converges {\tt
absolutely} and {\tt uniformly on any bounded subsets of
$\mathbb{R}\times C$}.
}%

\medskip

Following notations of (\ref{sdd13-4a}), we conclude that (A3) means
that for any $(t,\chi)\in \mathbb{R}\times C$ one has
$\Phi_d(t,\chi)\equiv \int^0_{-r} \chi (\theta) \, d
g_d(\theta,t,\varphi) = \sum_k \chi(\theta_k) \cdot h_k (t,\varphi)=
\sum_k \chi(-\eta_k(t,\varphi)) \cdot h_k (t,\varphi) $. Here all
$\eta_k$ and $h_k$ are {\tt continuous} functions.

\medskip

The first result is (c.f. \cite[lemma 1]{Rezounenko-CR-2011})

\smallskip

{\bf Theorem~1}. {\it Assume ${G} : \mathbb{R}\times X \times X \to
X$ is a continuous 
mapping and $p: \mathbb{R}\times X\to X$ (see (\ref{sdd13-3})) is
Lipschitz  $(||p(t,u)-p(s,v)||\le L_p (|s-t|+||u-v||)$, satisfying $
||p(s,u)||\le C_1||u|| +C_2, \forall (s,u)\in \mathbb{R}\times X$
with $C_i\ge 0$.
Under assumptions (A1)- (A3), the nonlinear mapping $B :
\mathbb{R}\times C \to X$, defined by (\ref{sdd13-2}), is
continuous.}

\medskip

{\bf Remark}. {\it It is important that nonlinear map $B$ is not
Lipschitz in the presence of discrete SDDs. 
The last means that discrete delays may be present, but be constant
or time-dependent only (i.e. $g_d(\theta,t,\varphi)= \widehat
g_d(\theta,t)$).

}

\medskip

{\it Proof of theorem 1.} Since a composition of continuous mappings
is continuous, it is enough (see (\ref{sdd13-2})) to show the
continuity of $F$ defined by (\ref{sdd13-3}).

We first split our $g$ in continuous and discontinuous parts
$g_c\equiv g_{ac}+g_s$ and $g_d$, respectively (see
(\ref{sdd13-4a})). This splitting gives the corresponding splitting
 $F=F_c+F_d,$ where $F_c$ corresponds to the continuous part
$g_c\equiv g_{ac}+g_s$.

{\it Case 1.} Let us first consider the {\it part $F_c$}.
We write
\begin{equation}\label{sdd13-t1-4}
F_c(t^1,\varphi) - F_c(t^2,\psi) = I_1 + I_2,
\end{equation}
where we denote
\begin{equation}\label{sdd13-t1-5}
I_1 =I_1(\varphi,\psi)\equiv \int^0_{-r}  \left[ p(t^1,\varphi
(\theta))-p(t^2,\psi(\theta))\right]  \, dg_c (\theta,t^1,\varphi),
\end{equation}

\begin{equation}\label{sdd13-t1-6}
I_2  =I_2(\varphi,\psi)\equiv \int^0_{-r} p(t^2,\psi(\theta)) \, d\,
[  g_c (\theta,t^1,\varphi) - g_c (\theta,t^2,\psi)].
\end{equation}
Using the Lipschitz property of $p$ and (A1), one can check that
\begin{equation}\label{sdd13-t1-7}
||I_1|| \le L_p (|t^1-t^2| + ||\varphi - \psi ||_C) \cdot 
 M_{Vg}.
\end{equation}
This  shows that $||I_1||\to 0$ when $|t^1-t^2| +||\varphi -\psi
||_C \to 0.$ To show that $||I_2||\to 0$ (when $t^1\to t^2$ and
$\varphi \to \psi$ in $C$) we use assumptions (A1) and (A2) to apply
the first Helly's theorem \cite[page 359]{Kolmogorov-Fomin}.

\smallskip

{\it Case 2.}  Now we prove the continuity of $F_d$ ({\it discrete
delays}).
Let us fix any $\varphi\in C, t\in \mathbb{R}$ and consider any
sequences $\{ \varphi^n \}\subset C$ and $\{ t^n \}\subset
\mathbb{R}$ such that $||\varphi^n-\varphi ||_C \to 0$ and $t^n\to
t$ when $n\to \infty$. Our goal is to prove that
$||F_d(t^n\varphi^n)-F_d(t,\varphi) || \to 0$.

Following the notations of (A3) 
we write
$$
F_d(t,\varphi) = \sum_k p(t,\varphi (-\eta_k(t,\varphi)))\cdot
h_k(t,\varphi)
$$
 and remind that it could be a series or a finite sum.
We split as follows
\begin{equation}\label{sdd13-t1-3}
F_d(t^n,\varphi^n)-F_d(t,\varphi) \equiv K^n_1 + K^n_2 +K^n_3 \in X,
\end{equation}
 where
$$K^n_1\equiv \sum_k  \left\{ p(t^n,\varphi^n
(-\eta_k(t^n,\varphi^n))) - p(t,\varphi
(-\eta_k(t^n,\varphi^n)))\right\} \cdot h_k(t^n,\varphi^n),$$
$$K^n_2\equiv \sum_k p(t,\varphi (-\eta_k(t^n,\varphi^n))) \cdot \left[ h_k(t^n,\varphi^n)-h_k(t,\varphi)\right],$$
$$K^n_3\equiv \sum_k  \left\{ p(t,\varphi
(-\eta_k(t^n,\varphi^n))) - p(t,\varphi
(-\eta_k(t,\varphi)))\right\} \cdot h_k(t,\varphi).$$
Using the Lipschitz property of $p$ one may check that
\begin{equation}\label{sdd13-t1-12}
||K^n_1|| \le L_p (|t^n-t|+ ||\varphi^n-\varphi ||_C ) \cdot \sum_k
|h_k(t^n,\varphi^n)|.
 \end{equation}
Now we discuss $K^n_2.$ The growth condition of $p$ implies
$||p(t,\varphi (-\eta_k(t^n,\varphi^n))) || \le (C_1 ||\varphi||_C +
C_2)$. Hence
 \begin{equation}\label{sdd13-t1-13}
 || K^n_2|| \le (C_1 ||\varphi||_C + C_2)\cdot \sum_k
 |h_k(t^n,\varphi^n)-h_k(t,\varphi)|.
 \end{equation}
In a similar way we obtain
\begin{equation}\label{sdd13-t1-14}
 || K^n_3|| \le L_p \sum_k
 |h_k(t,\varphi)|\cdot ||\varphi
(-\eta_k(t^n,\varphi^n)) - \varphi (-\eta_k(t,\varphi)||.
 \end{equation}
Now we show that $||K^n_j||\to 0$ as $n\to\infty$ for $j=1,2,3.$ The
first property $||K^n_1||\to 0$ follows from (A3) and
(\ref{sdd13-t1-12}). In (\ref{sdd13-t1-13}), the series converges
uniformly with respect to $n$ since the condition
$||\varphi^n-\varphi ||_C +|t^n-t|\to 0$  implies that $\{
(t,\varphi), (t^n,\varphi^n) \}$ is a bounded subset of
$\mathbb{R}\times C$. Assumption (A3) guarantees that each
$|h_k(t^n,\varphi^n)-h_k(t,\varphi)|$ is continuous with respect to
$(t^n,\varphi^n)$ and tends to zero when $n\to\infty$. Due to the
uniform convergence of the series in (\ref{sdd13-t1-13}) (see (A3)),
we arrive at $||K^n_2||\to 0$. To show that $||K^n_3||\to 0$ we also
mention that each $|h_k(t,\varphi)|\cdot ||\varphi
(-\eta_k(t^n,\varphi^n)) - \varphi (-\eta_k(t,\varphi)||$ (see
(\ref{sdd13-t1-14})) is continuous with respect to $(t^n,\varphi^n)$
and tends to zero as $n\to\infty$ due to (A3) and the strong
continuity of $\varphi\in C$. The uniform convergence (w.r.t.
$(t^n,\varphi^n)$) of the series in (\ref{sdd13-t1-14}) follows from
the estimate $|h_k(t,\varphi)|\cdot ||\varphi
(-\eta_k(t^n,\varphi^n)) - \varphi (-\eta_k(t,\varphi)||\le
|h_k(t,\varphi)|\cdot 2 ||\varphi ||_C$
(the right-hand side is independent of $n$!) and the Weierstrass 
dominant (uniform) convergence theorem. We conclude that
$||K^n_3||\to 0$. Since all $||K^n_j||\to 0$ as $n\to\infty$ for
$j=1,2,3$ we proved the property
$||F_d(t^n,\varphi^n)-F_d(t,\varphi) || \to 0$.
We shown that both $F_c$ and $F_d$ are continuous. The proof of
theorem~1 is complete. \finproof


\medskip

In our study we use the standard

\smallskip

{\bf Definition~1}. 
 {\it A function $u\in C([a-r,T]; X)$ is called a {\tt mild solution}
 on $[a-r,T]$ of the initial value problem (\ref{sdd13-1}), (\ref{sdd13-ic}) if it satisfies
 (\ref{sdd13-ic}) and
 \begin{equation}\label{sdd13-3-1}
u(t)=e^{-A (t-a)}\varphi(0) + \int^{t}_a e^{- A (t-s)} B (s,u_s)  \,
ds, \quad t\in [a,T].
 \end{equation}
}%

\medskip
\medskip

{\bf Theorem~2}. {\it Under the assumptions of theorem~1, the
initial value problem (\ref{sdd13-1}), (\ref{sdd13-ic})  possesses a 
mild solution for any $\varphi\in C$.}

\smallskip

The existence of a mild solution is a consequence of the continuity
of $B : \mathbb{R}\times C \to X$, given by theorem~1, which gives
us the possibility to use the standard method based on the Schauder
fixed point theorem (see \cite[theorem 3.1, p.4]{Fitzgibbon-JDE-1978}).


\medskip

{\bf Theorem~3}. {\it Let all the assumptions of theorem~1 are
valid. If additionally $||G(t,u,v)|| \le k_1(t)(||u||+||v||) +
k_2(t)$ with $k_i$ are locally integrable on $[a,\infty)$, then a
mild solution is global i.e. defined for all $t\ge a.$ }

\smallskip

The statement follows from theorem 2 and 
\cite[theorem 2.3, p. 49]{Wu_book}.

\bigskip

To get the uniqueness of mild solutions we need the following
additional assumptions.

\medskip

{\bf (A4)} {\it The total variation of function $g_c\equiv
g_{ac}+g_{s}$ satisfies 
\begin{equation}\label{sdd13-lip}
\exists L_{Vg_c}\ge 0 : \forall t^1,t^2\ge a \Rightarrow V^0_{-r}
[g_{c}(\cdot, t^1,\varphi) - g_{c}(\cdot, t^2,\psi)] \le L_{Vg_c}
 (|t^1-t^2|+||\varphi-\psi||_C).
\end{equation}
}%

\medskip

{\bf (A5)} {\it The discrete generating function $g_{d}$ satisfies
the following { uniform} condition:

\begin{itemize}
 \item there exists continuous function $\eta_{ign}(t)>0$, such that
 {\tt all} $\eta_k$ and $h_k$ "{\tt ignore}" values of
 $\varphi(\theta)$ for $\theta\in (-\eta_{ign}(t),0]$ i.e. 
 $$\hskip-12mm \exists\, \eta_{ign}(t)>0 : \forall t\ge a, 
 \forall\varphi^1, \varphi^2\in C :
 \forall\theta\in
 [-r,-\eta_{ign}(t)],\,\Rightarrow \varphi^1(\theta)= \varphi^2(\theta)\quad
   \Longrightarrow
   $$
 $$\forall k\in \mathbb{N} \Rightarrow \quad \eta_k (t,\varphi^1)=\eta_k (t,\varphi^2) 
 \quad
 \mbox{ and } \quad h_k (t,\varphi^1)= h_k (t,\varphi^2). 
 $$
\end{itemize}

}%

\medskip

{\bf Remark}. {\it Assumption (A5) is the natural generalization to
the non-autonomous case of multiple discrete state-dependent delays
of the condition introduced in \cite{Rezounenko_NA-2009}. In
\cite{Rezounenko_NA-2009,Rezounenko-CR-2011} the function
$\eta_{ign}(t)$ was constant $\eta_{ign}(t)\equiv \eta_{ign}>0$. For
more details and examples see \cite{Rezounenko_NA-2009} and also
\cite{Rezounenko-CR-2011}.
}%

\medskip

{\bf Theorem~4}. {\it Assume  (A1)- (A5) are valid, $p$ is as in
theorem~1, mapping ${G} : \mathbb{R}\times X \times X \to X$ is
continuous and locally Lipschitz with respect to its second and
third coordinates i.e.
for any
$R>0$ there exists $L_{G,R}>0$ such that for all
$t>a, ||u^i||,||v^i||\le R$ one has 
\begin{equation}\label{sdd13-i-6}
||G(t,u^1,v^1)-G(t,u^2,v^2)|| \le L_{G,R} \left(
||u^1-u^2||+||v^1-v^2||\right).
\end{equation}
 Then initial
value problem (\ref{sdd13-1}), (\ref{sdd13-ic})  possesses a 
{\tt unique} mild solution on an interval of the form $[a,b)$ where
$a<b\le +\infty$ for any $\varphi\in C$. The solution is continuous
with respect to initial data i.e. $||\varphi^n - \varphi||_C\to~0$
implies  $||u^n_t - u_t||_C\to 0$ for any $t\in [a,b).$ Here $u^n$
is the unique solution of (\ref{sdd13-1}), (\ref{sdd13-ic}) with
initial function $\varphi^n$ instead of $ \varphi$.
}%

\medskip

{\it Proof of theorem 4.} For the simplicity, we first consider a
particular case when the generating function $g=g_c\equiv
g_{ac}+g_{s}$ i.e. $F=F_c$ does not contain the discrete delays.

Let $(t^1,\varphi), (t^2,\psi)$ belong to a bounded subset
${\mathcal B}\subset \mathbb{R}\times C$. We use the splitting
(\ref{sdd13-t1-4}). One can see that
(see (\ref{sdd13-t1-6}))   
\begin{equation}\label{sdd13-t3-8}
 ||I_2|| \le M_{\mathcal B} \cdot  V^0_{-r} [ g_c (t^1,\varphi) -g_c (t^2,\psi) ].
\end{equation}
Assumption (A4) and (\ref{sdd13-t1-7}) imply that $F_c$ is locally
Lipschitz i.e. for any $R>0$ there exists $L_{F_c,R}>0$ such that
for all $||\varphi||_C\le R$, $||\psi||_C\le R$, $a\le t^1,t^2\le
a+R$ one has
\begin{equation}\label{sdd13-t3-9}
||F_c(t^1,\phi)- F_c(t^2,\psi)|| \le L_{F_c,R} \left(
|t^1-t^2|+||\phi-\psi||_C \right).
\end{equation}


Consider a sequence $\{ \varphi^n \}\subset C$ such that
$||\varphi^n-\varphi ||_C \to 0$ as $n\to\infty$. Denote by
$u(t)=u(t;\varphi)$ 
{\it any} mild solution of (\ref{sdd13-1}), (\ref{sdd13-ic})  and by
$u^n(t)=u^n(t;\varphi^n)$ {\it any} mild solution of
(\ref{sdd13-1}), (\ref{sdd13-ic}) with initial data $\varphi^n\in
C$. The existence of these solutions is proved in theorem 2. The
Schauder fixed point theorem (see e.g. \cite[theorem 2.1,
p.46]{Wu_book}), used in the proof of theorem 2 implies that one can
choose $R>0$ and $T\in [a,b)$ to have
 $||\varphi^n||_C\le R$, $||u^n(t;\varphi^n)||\le R$ for all $t\in [a,T]$.

 Using the local Lipschitz property of mapping $G$ (\ref{sdd13-i-6}), (\ref{sdd13-t3-9}) %
 and  the form (\ref{sdd13-2}), we get
\begin{equation}\label{sdd13-i-7} ||B(t,\varphi)- B(t,\psi)|| \le L_{G,R}(1+L_{F_c,R})
||\varphi-\psi||_C  + L_{G,R}||F_d(t,\varphi)- F_d(t,\psi)||.
\end{equation}

Hence for any $t\in [a,T]$ one has (we remind that $||T(t)||\equiv
||e^{-At}||\le  e^{\omega t}$ and $F=F_c$)
$$ ||u_t-u^n_t||_C \le  e^{\omega (T-a)}||\varphi -\varphi^n||_C + L_{G,R}(1+L_{F_c,R})    e^{\omega (T-a)}\cdot
\int^t_a ||u_s-u^n_s||_C \, ds.
$$
The last estimate (by the Gronwall lemma) implies
$$ ||u_t-u^n_t||_C \le
 e^{L_{G,R}(1+L_{F_c,R})    e^{\omega (T-a)} }  e^{\omega (T-a)} \cdot ||\varphi -\varphi^n||_C.
$$
That is
\begin{equation}\label{sdd13-CT}
 ||u_t-u^n_t||_C \le
C_T\cdot ||\varphi -\varphi^n||_C, \quad \forall t\in [a,T],\, 
C_T \equiv
  e^{\omega (T-a)} \exp\{ L_{G,R}(1+L_{F_c,R})    e^{\omega (T-a)}
  \}.
\end{equation}
It proves the uniqueness of mild solutions and the continuity with
respect to initial data in the case $g=g_c$.

The second particular case $g=g_d$ (the purely discrete delay) and
only one point of discontinuity has been considered in detail in
\cite{Rezounenko_NA-2009} (the autonomous case). It was proved in
\cite{Rezounenko_NA-2009} that (A5)  implies the desired result.

Now we consider the general case (both discrete and distributed 
delays, including the case of multiple discrete SD-delays).

Using the splitting $F=F_d+F_c$, we have, by definition of mild solutions,
$$u^n(t)-u(t) = e^{-A(t-a)} (\varphi^n(0)-\varphi(0)) +
\int^t_a e^{-A(t-\tau)} \left\{F_d(\tau,u^n_\tau) -
F_d(\tau,u_\tau)\right\}\, d\tau $$ $$ + \int^t_a e^{-A(t-\tau)}
\left\{F_c(\tau,u^n_\tau) - F_c(\tau,u_\tau)\right\}\, d\tau.
$$
Using (\ref{sdd13-t3-9}), one gets for all $t\in [a,T]$
$$||u^n(t)-u(t)|| \le ||\varphi^n(0)-\varphi(0)||  e^{\omega (T-a)} +
 e^{\omega (T-a)} L_{G,R}\int^t_a ||F_d(\tau,u^n_\tau) -
F_d(\tau,u_\tau)||\, d\tau
$$ $$
+ L_{G,R}(1+L_{F_c,R})  e^{\omega (T-a)}\int^t_a ||u^n_\tau -
u_\tau||_C\, d\tau
$$
\begin{equation}\label{sdd13-t3-10}
= G^n(t) + L_{G,R}(1+L_{F_c,R})  e^{\omega (T-a)}\int^t_a ||u^n_\tau
- u_\tau||_C\, d\tau,
\end{equation}
where
\begin{equation}\label{sdd13-t3-12}
G^n(t)\equiv ||\varphi^n(0)-\varphi(0)||  e^{\omega (T-a)} + L_{G,R}
e^{\omega (T-a)} \int^t_a ||F_d(\tau,u^n_\tau) -
F_d(\tau,u_\tau)||\, d\tau
\end{equation}
 is a nondecreasing (in time) function.
Multiply the last estimate by $e^{-t L_{G,R}(1+L_{F_c,R})  e^{\omega
(T-a)} }$ to get
$$ {d\over dt} \left( e^{-t L_{G,R}(1+L_{F_c,R})  e^{\omega
(T-a)} } \int^t_a ||u^n_\tau - u_\tau||_C\, d\tau\right) \le e^{-t
L_{G,R}(1+L_{F_c,R})  e^{\omega (T-a)} } G^n(t),
$$
which, after integration from $a$ to $t$, shows that ($G^n(t)$ is
nondecreasing)
$$ e^{-t L_{G,R}(1+L_{F_c,R})  e^{\omega
(T-a)} } \int^t_a ||u^n_\tau - u_\tau||_C\, d\tau \le \int^t_a
e^{-\tau L_{G,R}(1+L_{F_c,R})  e^{\omega (T-a)} } G^n(\tau)\, d\tau
$$
$$\le G^n(t)\int^t_a e^{-\tau L_{G,R}(1+L_{F_c,R})  e^{\omega (T-a)} } \, d\tau
$$
$$= G^n(t) \left( e^{-a L_{G,R}(1+L_{F_c,R})  e^{\omega (T-a)} }-e^{-t
L_{G,R}(1+L_{F_c,R})  e^{\omega (T-a)} }\right)
[L_{G,R}(1+L_{F_c,R})  ]^{-1} e^{-\omega (T-a)}.
$$

We have
$$ L_{G,R}(1+L_{F_c,R})  e^{\omega (T-a)}\int^t_a ||u^n_\tau - u_\tau||_C\, d\tau
\le G^n(t) \left( e^{(t-a) L_{G,R}(1+L_{F_c,R})  e^{\omega (T-a)}
}-1 \right).
$$
We substitute the last estimate into (\ref{sdd13-t3-10}) to obtain
\begin{equation}\label{sdd13-t3-11}
||u^n_t-u_t||_C \le  G^n(t)\cdot e^{(t-a) L_{G,R}(1+L_{F_c,R})
e^{\omega (T-a)} }.
\end{equation}

Let us fix any $c >b$ and denote by $\eta_{ign}\equiv \min \{
\eta_{ign}(t) : t\in [a,c]\}$. By Assumption (A5), $\eta_{ign}(t)>0$
and continuous, so $\eta_{ign}>0$. Let us denote by $\sigma\equiv
\min\{b,a+\eta_{ign}\}>a.$

 Now our goal is to show that for any fixed $t\in [a, \sigma)$
one has $G^n(t)\to 0$ when $n\to\infty$ (we remind that
$||\varphi^n-\varphi ||_C \to 0$).

Let us consider the extension functions
$$
\overline \varphi(s)\equiv \left[\begin{array}{ll}
  \varphi(s) & s \in [-r, 0]; \\
  \varphi(0) & s\in (0, \sigma) \\
\end{array}
\right. \quad \hbox{ and }\quad  \overline \varphi^n(s)\equiv
\left[\begin{array}{ll}
  \varphi^n(s) & s \in [-r, 0]; \\
  \varphi^n(0) & s\in (0, \sigma) \\
\end{array}.
\right.
$$

An important consequence of (A5) is that
$F_d(t,u_t)=F_d(t,\bar{\varphi}_{t})$ for all $t\in [a,\sigma)$ and
any solution $u:[a-r,\sigma)\to X$, satisfying $u_a=\varphi$. In the
same way $F_d(t,u^n_t)=F_d(t,\bar{\varphi}^n_{t})$ for all $t\in
[a,\sigma)$. Hence the continuity of $F_d$ implies
$||F_d(\tau,\overline \varphi^n_\tau) -F_d(\tau,\overline
\varphi_\tau)||\to 0$ for any $\tau\in [a, \sigma)$.

\medskip

{\bf Remark}. {\it We notice that the case we consider now is
simpler that the one in the proof of theorem 1 (see
(\ref{sdd13-t1-3})) since we estimate $F_d$ at the same first
coordinate (time moment $\tau\in [a, \sigma)$).

}%

\medskip

 The property $||F_d(\tau,\overline \varphi^n_\tau) -F_d(\tau,\overline
\varphi_\tau)||\to 0$ for any $\tau\in [a, \sigma)$ and the uniform
boundedness of the term allows us to use the classical
Lebesgue-Fatou lemma (see \cite[p.32]{yosida}) for the scalar
function $||F_d(\tau,\overline \varphi^n_\tau) -F_d(\tau,\overline
\varphi_\tau)||$ to conclude that $G^n(t)\to 0$ when $n\to\infty$
(for any fixed $t\in [a, \sigma)$). Hence (\ref{sdd13-t3-11}) gives
the continuity of the mild solutions with respect to initial
functions for all $t\in [a, \sigma)$. Particularly, it gives the
uniqueness of solutions. For bigger time values we use the chain
rule (by the uniqueness) for steps less than or equal to, say
$(\sigma-a)/2$ (for more details see \cite{Rezounenko_NA-2009}).
Since the composition of continuous mappings is continuous, the
proof of theorem~4 is complete. \finproof

\medskip

{\bf Remark}. {\it Discussing the proof of theorem 4, we see that in
the case $g=g_c$ (no discrete state-dependent delays) the delay
mapping $F=F_c$ is locally Lipschitz continuous (see
(\ref{sdd13-t3-9})) and, as a consequence, one has a standard
estimate for the difference of two solutions (\ref{sdd13-CT}) in
terms of the difference of initial functions
$||\varphi-\varphi^n||_C$. In case of the presence of discrete SDDs
we have estimate (\ref{sdd13-t3-11}) with $G^n$ defined by
(\ref{sdd13-t3-12}) and for the Lebesgue-Fatou lemma it was enough
to have property $||F_d(\tau,\overline \varphi^n_\tau)
-F_d(\tau,\overline \varphi_\tau)||\to 0$ which does not provide
information on the difference $||F(t,\varphi^n)-F(t,\varphi)||$ in
terms of $||\varphi-\varphi^n||_C$ (it is definitely not a Lipschitz
property). A way to get such an information is to use the modulus of
continuity $\omega_f(\delta;Y)$. We remind that
$\omega_f(\delta;Y)\equiv \sup \{ ||f(x)-f(y)|| : x,y\in Y,
||x-y||\le \delta \}$. For the simplicity of presentation we
consider $F_d$ with one discrete SDD. We have (see
(\ref{sdd13-t1-3}) with $t=t^n=\tau$)
\begin{equation}\label{sdd13-t3-13}
F_d(\tau,\varphi^n) -F_d(\tau,\varphi)= p(\tau,\varphi^n
(-\eta(\tau,\varphi^n)))\cdot h(\tau,\varphi^n) -p(\tau,\varphi
(-\eta(\tau,\varphi)))\cdot h(\tau,\varphi).
\end{equation}
The estimates for $K^n_i, i=1,2,3$ (see (\ref{sdd13-t1-12})-
(\ref{sdd13-t1-14}) with $k=1$) show that
$$||F_d(\tau,\varphi^n) -F_d(\tau,\varphi)|| \le  L_p V^0_{-r}g_d
\cdot ||\varphi-\varphi^n||_C  $$ $$+ (C_1||\varphi||_C+C_2)\cdot
\omega_h\left(||\varphi-\varphi^n||_C;\tilde{Y}\right) + L_p
V^0_{-r}g_d \cdot \omega_\varphi \left( \omega_\eta
\left(||\varphi-\varphi^n||_C;\tilde{Y}\right); [-r,0]\right),
$$
where we denoted by $\tilde{Y}\equiv \{ (t,\bar{\varphi}_t),
(t,\bar{\varphi}^n_t) : t\in [a,\sigma_1], n\in \mathbb{N} \}$, with
$a<\sigma_1<\sigma$.  By (A1) one has $V^0_{-r}g_d \le M_{Vg}$ and
\begin{equation}\label{sdd13-t3-14}
||F_d(\tau,\varphi^n)
-F_d(\tau,\varphi)|| \le  L_p M_{Vg} \cdot
\left[||\varphi-\varphi^n||_C  + \omega_\varphi \left( \omega_\eta
\left(||\varphi-\varphi^n||_C;\tilde{Y}\right);
[-r,0]\right)\right]$$ $$+ (C_1||\varphi||_C+C_2)\cdot
\omega_h\left(||\varphi-\varphi^n||_C;\tilde{Y}\right).
\end{equation}
Since $\varphi^n\to \varphi$ in $C$, we see that $\tilde{Y}$ is
compact. Using (A3) and the classical Cantor theorem, we know that
$h$ and $\eta$ are equicontinuous on $\tilde{Y}$ and
$\omega_h\left(||\varphi-\varphi^n||_C;\tilde{Y}\right)\to 0$ and
$\omega_\eta\left(||\varphi-\varphi^n||_C;\tilde{Y}\right)\to 0$ as
$||\varphi-\varphi^n||_C \to 0$. We remind that
$F_d(t,u_t)=F_d(t,\bar{\varphi}_{t})$,
$F_d(t,u^n_t)=F_d(t,\bar{\varphi}^n_{t})$ and use
$||\bar{\varphi}_t-\bar{\varphi}^n_t||_C\le
||\varphi-\varphi^n||_C$. Finally,  one can substitute the estimate
(\ref{sdd13-t3-14}) into (\ref{sdd13-t3-12}) and then the estimate
for $G^n$ into (\ref{sdd13-t3-11}) to get an estimate for the
difference of two solutions in terms of the difference of initial
functions $||\varphi-\varphi^n||_C$.
}%

\medskip

{\bf Remark}. {\it In the theory of {\tt ordinary} differential equations with SDDs
it is usual  to restrict the class of initial functions $\varphi$ to Lipschitz ones
\cite{Walther_JDE-2003,Hartung-Krisztin-Walther-Wu-2006}. In this case one 
restricts the set of SDDs (both $\eta$ and $h$) to Lipschitz mappings. In such a situation the previous remark
evidently provides the Lipschitz property of $F_d$ and hence $F$. It follows 
 from the property
$f\in {\mathcal Lip}(L_f;Y) \Longrightarrow \omega_f(\delta;Y)\le L_f \cdot\delta$ and the estimates above.

}%

\section{Invariance}

This section is devoted to an extension of the fundamental invariance principle \cite{Martin-Smith-TAMS-1990}
to the case of PDEs with discrete state-dependent delays.

We also refer the reader to \cite{Ruess-TAMS-2009} for important
generalizations of the invariance principle in several directions
(without SDDs) and for a survey of the existing literature on the
subject. In the present paper the emphasis is on the delayed term,
not on the partial differential operator.

Following \cite{Martin-Smith-TAMS-1990}, we assume that the next
hypotheses are satisfied:

\begin{description}
  \item[(H1)] $D$ is a closed subset of $[a-r,\infty)\times X$ and $D(t)\equiv \{ x\in X : (t,x)\in D\}$ is nonempty for each $t\ge a-r.$
  \item[(H2)] ${\mathcal D}$ is the closed subset of $[a,\infty)\times C$ defined by ${\mathcal D}\equiv \{ (t,\varphi) : {\varphi (\theta)\in D(t+\theta)} \hbox{ for all } {-r\le\theta\le 0} \}.$  Also, ${\mathcal D}(t)\equiv \{ \varphi \in C : (t,\varphi) \in {\mathcal D}\}$ for each $t\ge a,$ and we assume that ${\mathcal D}(t)$ is nonempty for each set $t\ge a.$
  \item[(H3)] For each $b>a$ there are a $\hat K(b)>0$ and a continuous nondecreasing function $\eta_b : [0,b-a)\to [0,\infty)$ satisfying $\eta_b(0)=0$ with the property that if $a\le t_1<t_2\le b, x_1\in D(t_1),$ and $x_2\in D(t_2),$ then there is a continuous function $w:[t_1,t_2]\to X$ such that $w(t_1)=x_1, w(t_2)=x_2, w(t)\in D(t)$ for $t_1<t<t_2,$ and
      $$ |w(t)-w(s)| \le \eta_b(|t-s|) + \hat K(b) |t-s| \frac{|x_2-x_1|}{t_2-t_1}
      $$
      for all $s,t\in [t_1,t_2].$
  \item[(H4)] $ B$ is continuous from $D(B)$ into $X$ where ${\mathcal D}\subset D(B)\subset [a,\infty)\times C.$
\end{description}

\medskip

{\bf Remark} \cite[page 16]{Martin-Smith-TAMS-1990}. {\it If $D$ is
{\tt convex} then (H3) is automatically satisfied by defining
$$w(t)=\frac{(t_2-t)x_1+(t-t_1)x_2}{(t_2-t_1)}\quad \hbox{ for } \quad  t_1\le t\le
t_2.$$ We see that $||w(t)-w(s)||=\left\|
\frac{(s-t)x_1+(t-s)x_2}{(t_2-t_1)} \right\| \le
\frac{||x_1-x_2||}{(t_2-t_1)} |t-s|$ for $t_1\le s<t\le t_2$ and
hence (H3) is satisfied with $\hat K(b)\equiv 1$ and $\eta_b=0.$
}%

\medskip

We use the notation
$$ d(x;D(t))\equiv \inf \{ |x-y| : y\in D(t)\} \quad \hbox{for } \quad  x\in X, t\ge a.
$$
The fundamental criterion for the invariance of the set ${\mathcal
D}$, called the {\it subtangential condition}, (see
\cite[(2.2)]{Martin-Smith-TAMS-1990}) is given by
\begin{equation}\label{sdd13-i-1}
\lim\limits_{h\to 0+} \frac{1}{h}\, d\left( e^{-Ah}\varphi(0) + \int^{t+h}_t  e^{-A(t+h-s)} B(t,\varphi)\, ds\, ; D(t+h)  \right)\quad \hbox{for }\quad  (t,\varphi)\in {\mathcal D}.
\end{equation}

\medskip

The following result is an extension of \cite[theorem
2]{Martin-Smith-TAMS-1990} to the case of general state-dependent
delay.

\medskip

{\bf Theorem~5}. {\it Assume  (A1)- (A5), (H1)-(H4) and
(\ref{sdd13-i-1}) are valid. Let mapping   $p$ be as in theorem~1
and $G  : \mathbb{R}\times X \times X \to X$ satisfy the property:
for each $R>0$ there are 
an $L_{G,R}>0$ and a continuous $\nu_R:[0,\infty)\to [0,\infty)$ 
such that $\nu_R(0)=0$ and
\begin{equation}\label{sdd13-i-12}
||G(t^1,u^1,v^1)-G(t^2,u^2,v^2)|| \le \nu_R(|t^1-t^2|) + L_{G,R}
\left( ||u^1-u^2||+||v^1-v^2||\right)
\end{equation}
for all $||u^i||, ||v^i||\le R$ and $a\le t^1, t^2\le a+R$.

Then initial value problem (\ref{sdd13-1}), (\ref{sdd13-ic})  has a 
unique mild solution on an interval of the form $[a,b)$ where $a<b\le +\infty$ for any $\varphi\in C.$
If additionally $\varphi\in {\mathcal D}(a)$, then $u(t)\in D(t)$ for
$a\le t<b$ and if $b<+\infty$ then $||u_t||_C\to \infty$ as $t\to b-0$.
}%

\medskip

{\bf Remark}. {\it In contrast to \cite[theorem
2]{Martin-Smith-TAMS-1990}, the nonlinear term $B$ in equation
(\ref{sdd13-1}) is not Lipschitz (with respect to the second
coordinate, c.f. \cite[property (2.3),
p.18]{Martin-Smith-TAMS-1990}) in the presence of discrete SDD.
So \cite[theorem 2]{Martin-Smith-TAMS-1990} could not be applied to
our case. Instead, assumption (A5) provides the uniqueness of mild
solutions and saves the line of the proof presented in
\cite{Martin-Smith-TAMS-1990}.
}%

\medskip

{\it The proof of theorem~5} follows closely that of \cite[theorem
2]{Martin-Smith-TAMS-1990}. The last
consists of eight lemmas and the final
part which spends pages 35-43 of the original article. There is no
need to repeat these lemmas since they are not affected by the lack
of the Lipschitz property of $B$ and we refer the reader to
\cite{Martin-Smith-TAMS-1990} for all notations and details.
 Here we only remind the main steps of the original proof and give
 the new part of the proof based on assumption (A5).

First, for fixed $\sigma>a, \varepsilon_0>0$ and any $\varepsilon\in
[0,\varepsilon_0]$ the $\varepsilon$-approximate solution $w$ is
constructed. It is done by a careful construction (see
\cite{Martin-Smith-TAMS-1990} for all details) of an increasing
sequence $\{ t_i \}^\infty_0 \subset [a,\sigma+ \varepsilon_0]$ such
that $w(t_i)\in D(t_i)$ and (see
\cite[(4.6)]{Martin-Smith-TAMS-1990})
\begin{equation}\label{sdd13-i-2}
\left\| e^{-A(t_{i+1}-t_i)}w(t_{t_i}) + \int^{t_{i+1}}_{t_i}
e^{-A(t_{i+1}-s)} B(t_i,w_{t_i})\, ds\, - w(w_{t_{i+1}}) \right\|
\quad \le \varepsilon (t_{i+1}-t_{i}).
\end{equation}

Let $\{ \varepsilon_n\}^\infty_{1}$ be a decreasing sequence such
that $\varepsilon_n\to 0$ as $n\to\infty$ and for each $n\ge 1$ let
$w^n$ and $\{ t^n_i \}^\infty_{i=0} $ be as constructed above with
$\varepsilon = \varepsilon_n, t_i=t^n_i,$ and $w=w^n$. Denote by
$\gamma^n : [a,\sigma]\to [a,\sigma]$ the function $\gamma^n
(t)=t^n_i$ whenever $t\in [t^n_i,t^n_{i+1}]$.

\medskip

{\bf Remark~A}. {\it In addition to consideration in
\cite{Martin-Smith-TAMS-1990}, we assume that $\sigma+ \varepsilon_0
-a< \eta_{ign}\equiv\min_{t\in[a,c]}\eta (t)$ for some fixed $c>a$.
Since we prove the local existence, $c$ could be chosen arbitrary
and (A5) gives $\eta_{ign}>0$ ($\eta (t)$ is continuous).
}%

\medskip

For convenience a companion function $v^n$ for $w^n$ is defined in
the following manner (see \cite[(4.9)]{Martin-Smith-TAMS-1990})
\begin{equation}\label{sdd13-i-3}
v^n(t)= e^{-A(t-a)}\varphi(0) + \int^{t}_{a} e^{-A(t-s)}
B(\gamma^n(s),w^n_{\gamma^n(s)})\, ds\, \quad \hbox{for } t\in
[a,\sigma].
\end{equation}
and $v^n(a+\theta)=\varphi(\theta)$ for $\theta\in [-r,0]$.

It is shown (see \cite[lemmas 4.6 and 4.7)]{Martin-Smith-TAMS-1990})
that
\begin{equation}\label{sdd13-i-4}
||v^n(t)-w^n(t)|| \le \hat{P} \max \{
\varepsilon_n,\eta_b(\varepsilon_n)\} \quad \hbox{for}\quad  t\in
[a-r,\sigma],\, n=1,2,\ldots, \, \hat{P}>0,
\end{equation}
\begin{equation}\label{sdd13-i-5}
||w^n_t-w^n_{\gamma^n(t)}||_C \le \hat{Q} \max \{
\varepsilon_n,\eta_b(\varepsilon_n)\} \quad \hbox{for}\quad  t\in
[a,\sigma],\, n=1,2,\ldots, \, \hat{Q}>0,
\end{equation}
with $\hat{P},\hat{Q}>0$ both independent of $t$ and $n$.

Next, \cite[lemma 4.8]{Martin-Smith-TAMS-1990} shows that {\it if}
there is a function $u:[a-r,\sigma]\to X$ such that
$u(t)=\lim_{n\to\infty} w^n(t)$ uniformly for $t\in [a-r,\sigma]$,
{\it then} $(t,u(t))\in D$ and $u$ is a {\it solution} to
(\ref{sdd13-1}), (\ref{sdd13-ic}) on $[a,\sigma]$.

\medskip

Now we proceed the final part of the proof where we use assumption
(A5) instead of the Lipschitz property of $B$. First, we need an
estimate for $||B(t^1,u_{t^1})- B(t^2,u_{t^2})||$. We consider the
splitting of the delay mapping $F$ onto continuous and discrete
parts $F=F_c+F_d$ (according to (\ref{sdd13-4a})) and use the local
Lipschitz property of $F_c$ (due to (A4), see (\ref{sdd13-t3-9})).
 More precisely, using the local Lipschitz property of mapping $G$ (\ref{sdd13-i-12}), (\ref{sdd13-t3-9}) %
 and  the form (\ref{sdd13-2}), we get
$$||B(t^1,u_{t^1})- B(t^2,u_{t^2})||
$$
$$\le \nu_R\left(
|t^1-t^2|\right) + L_{G,R} L_{F_c,R}\cdot |t^1-t^2| +
L_{G,R}(1+L_{F_c,R})\cdot  ||u_{t^1}-u_{t^2}||_C
$$
\begin{equation}\label{sdd13-i-11}
+ L_{G,R}\, ||F_d(t^1,u_{t^1})- F_d(t^2,u_{t^2})||.
\end{equation}
Since $|t-\gamma^n(t)|\le \varepsilon_n$ it follows that
\begin{equation}\label{sdd13-i-8}
|\gamma^n(t)-\gamma^m(t)|\to 0 \quad \hbox{ as } \quad n,m\to\infty
\quad \hbox{ uniformly for }\quad t\in [a,\sigma].
\end{equation}
Using (\ref{sdd13-i-4}), (\ref{sdd13-i-5}), consider $\bar{R}>0$
such that $||w^n(t)||\le \bar{R}$ for all $n\ge 1$ and $t\in
[a-r,\sigma]$. Then (\ref{sdd13-i-11}) implies



$$||B(\gamma^n(s),w^n_{\gamma^n(s)})- B(\gamma^m(s),w^m_{\gamma^m(s)})||
$$
$$\le \nu_R\left(
|\gamma^n(s)-\gamma^m(s)|\right) + L_{G,R} L_{F_c,R}\cdot
|\gamma^n(s)-\gamma^m(s)| + L_{G,\bar{R}}(1+L_{F_c,R})
||w^n_{\gamma^n(s)}-w^m_{\gamma^m(s)}||_C
$$
    $$+ L_{G,\bar{R}}||F_d(\gamma^n(s),w^n_{\gamma^n(s)})-
F_d(\gamma^m(s),w^m_{\gamma^m(s)})|| \le [\hbox{estimates
(\ref{sdd13-i-4}), (\ref{sdd13-i-5}) give} ]
$$
$$ \le L_{G,\bar{R}}(1+L_{F_c,R}) ||v^n_s-v^m_s ||_C + L_{G,\bar{R}}||F_d(\gamma^n(s),w^n_{\gamma^n(s)})-
F_d(\gamma^m(s),w^m_{\gamma^m(s)})|| +
\widetilde{\varepsilon}_{n,m},
$$ where $\widetilde{\varepsilon}_{n,m}\to 0$ as $n,m\to\infty$ (due to
(\ref{sdd13-i-8})). Using (\ref{sdd13-i-3}), we obtain
$$ ||v^n(t)-v^m(t)|| \le \int^t_a M L_{G,\bar{R}}(1+L_{F_c,R})||v^n_s-v^m_s
||_C\, ds $$
\begin{equation}\label{sdd13-i-9}
+ \int^t_a  M L_{G,\bar{R}}||F_d(\gamma^n(s),w^n_{\gamma^n(s)})-
F_d(\gamma^m(s),w^m_{\gamma^m(s)})||\, ds + \hat{\varepsilon}_{n,m},
\end{equation}
where $\hat{\varepsilon}_{n,m}\to 0$ as $n,m\to\infty$.

Now our goal is to show that the second integral in
(\ref{sdd13-i-9}) tends to zero as $n,m\to\infty$. Let us denote by
$\bar{\varphi}$ the extension function
$\bar{\varphi}(a+\theta)=\varphi(\theta)$ for $\theta\in [-r,0]$ and
$\bar{\varphi}(s)=\varphi(0)$ for $s\in (a,\sigma]$. We remind that
$\sigma -a< \eta_{ign}$ (see remark~A above). An important
consequence of (A5) is that $F_d(t,u_t)=F_d(t,\bar{\varphi}_{t})$
for all $t\in [a,\sigma]$ and any continuous function
$u:[a-r,\sigma]\to X$, satisfying $u_a=\varphi$. Since all $w^n$, by
construction, are continuous and satisfy $w^n_a=\varphi$, we arrive
to the property $F_d(t,w^n_t)=F_d(t,\bar{\varphi}_{t})$ for all
$t\in [a,\sigma].$ Hence (see the second integral in
(\ref{sdd13-i-9}))
$$G^{n,m}_d(s)\equiv ||F_d(\gamma^n(s),w^n_{\gamma^n(s)})-
F_d(\gamma^m(s),w^m_{\gamma^m(s)})|| =
||F_d(\gamma^n(s),\bar{\varphi}_{\gamma^n(s)})-
F_d(\gamma^m(s),\bar{\varphi}_{\gamma^m(s)})||.
$$
The continuity of $F_d$ and (\ref{sdd13-i-8}) give $G^{n,m}_d(s)\to
0$ as $n,m\to\infty$ for all $s\in [a,\sigma].$ Since $G^{n,m}_d(s)$
is bounded, the classical Lebesgue-Fatou lemma implies that
$$\int^t_a G^{n,m}_d(s)\, ds \to 0 \hbox{ as } n,m\to\infty.
$$
The last property gives (see (\ref{sdd13-i-9}))
\begin{equation}\label{sdd13-i-10}
||v^n(t)-v^m(t)|| \le \int^t_a M
L_{G,\bar{R}}(1+L_{F_c,R})||v^n_s-v^m_s ||_C\, ds +
\varepsilon_{n,m},
\end{equation}
where $\varepsilon_{n,m}\to 0$ as $n,m\to\infty$.

The rest of the proof follows \cite[page
43]{Martin-Smith-TAMS-1990}. Defining $q_{n,m}(t)\equiv \max \{
||v^n(s)-v^m(s)|| : a-r\le s\le t\}$ we see that for each $t\in
[a,\sigma]$ there is an $\alpha(t)\in [a-r,t]$ such that
$$q_{n,m}(t) = ||v^n(\alpha(t))-v^m(\alpha(t))||
\le \int^{\alpha(t)}_a M L_{G,\bar{R}}(1+L_{F_c,R})||v^n_s-v^m_s
||_C\, ds + \varepsilon_{n,m}
$$
$$ \le \int^{\alpha(t)}_a M L_{G,\bar{R}}(1+L_{F_c,R})\, q_{n,m}(s) \, ds +
\varepsilon_{n,m}.
$$
Gronwall's inequality along with the fact that $\varepsilon_{n,m}\to
0$ as $n,m\to\infty$ shows that $q_{n,m}(t)\to 0$ as $n,m\to\infty$,
and hence $\{ v^n(t)\}^\infty_{n=1}$ is uniformly Cauchy on
$[a-r,\sigma]$. This implies that $\{ w^n(t)\}^\infty_{n=1}$ is
uniformly Cauchy on $[a-r,\sigma]$ and hence initial value problem
(\ref{sdd13-1}), (\ref{sdd13-ic})  has a 
mild solution on $[a-r,\sigma]$ (see the discussion above and
\cite[lemma 4.8]{Martin-Smith-TAMS-1990}). The uniqueness of
solution is due to (A5) and provided by theorem 3. The standard
continuation arguments give solutions defined on a maximal interval.
The proof of theorem 5 is complete. \finproof

\medskip

The following important corollary remains valid in the presence of
SDDs.

\medskip

{\bf Corollary} (c.f. \cite[page 18]{Martin-Smith-TAMS-1990}). {\it
Suppose $K$ is a closed, convex subset of $X$ and all the
assumptions of theorem 5 are satisfied with $D(t)\equiv K$ for all
$t\ge a.$ Suppose further that
\begin{description}
    \item[(a)] $T(t): K\to K$ for $t\ge 0$ and
    \item[(b)] $\lim_{h\to 0+} {1\over h} d(\varphi(0)+h B(t,\varphi);
    K) =0$ for $(t,\varphi)\in {\mathcal D}$.
\end{description}

Then (\ref{sdd13-1}), (\ref{sdd13-ic})  has a unique noncontinuable
mild solution $u$ on $[a,b)$ for some $b>a$ and $u(t)\in K$ for all
$t\in [a-r,b)$.
}%

\medskip

Since we are interested in models from biology (see the examples
above) the following remark is of prime importance for us.

\medskip

{\bf Remark} (c.f. \cite[page 7]{Martin-Smith-TAMS-1990}). {\it
Consider system (\ref{sdd13-ex-2-a}). If $K=[0,\infty)^m$, then
condition (a) of the previous corollary holds and condition (b)
holds only in case $G=(G^i)^m_1$ is quasipositive: if $k\in \{
1,\ldots,m\}$ and $(t,\varphi)\in [a,\infty)\times
C([-r,0];{C(\bar{\Omega})}^m)$ with $\varphi^i(\theta,x)\ge 0$ for
all $-r\le \theta\le 0, x\in \bar{\Omega}$ and $i=1,\ldots,m,$ then
$\varphi^i(0,\cdot)=0$ implies $G^i(t,F^i(t,\varphi(\cdot,x) )) \ge
0$ for all $x\in \bar{\Omega}$. This condition gives criteria to
determine if solutions to (\ref{sdd13-ex-2-a}) remain nonnegative if
they are nonnegative initially.
}%
$$
$$

\medskip



\medskip





\begin{thebibliography}{9}



\bibitem
{Babin-Vishik}
 A.V. Babin, and M.I.~Vishik, "Attractors of Evolutionary
 Equations", Amsterdam, North-Holland, 1992.

\bibitem
{NA-1998} L. Boutet de Monvel, I.D. Chueshov and A.V. Rezounenko,
Inertial manifolds for retarded semilinear parabolic equations,
Nonlinear Analysis, 34 (1998) 907-925.




\bibitem
{Chueshov-JSM-1992} I.~D. Chueshov, On a certain system of equations
with delay, occuring in aeroelasticity, J. Soviet Math. { 58},
(1992) 385-390.

\bibitem
{Cras-1995}
 I.~D. Chueshov and  A.~V. Rezounenko, Global attractors for a class of
retarded quasilinear partial differential equations,
C.R.Acad.Sci.Paris, Ser.I {\bf 321}, 607-612 (1995);  (~detailed
version: Math.Physics, Analysis, Geometry, Vol.2,  N.3 (1995),
363-383).

\bibitem
{Chueshov_book}
 I.~D. Chueshov, "Introduction to the Theory of
Infinite-Dimensional Dissipative Systems", Acta, Kharkov, 1999) (in
Russian). English transl. Acta, Kharkov (2002) (see {\small
 http://www.emis.de/monographs/Chueshov }).




\bibitem
{Walther_book} O. Diekmann, S. van Gils, S. Verduyn Lunel, H-O.
Walther,  "Delay Equations: Functional, Complex, and Nonlinear
Analysis",  Springer-Verlag, New York, 1995.

\bibitem%
{Fitzgibbon-JDE-1978}
W.E. Fitzgibbon,  Semilinear functional differential equations in Banach space,
 J. Differential Equations 29 (1978), no. 1, 1--14. 



\bibitem
{Hadamard-1902} J.~Hadamard, Sur les probl\`{e}mes aux derivees
partielles et leur signification physique, Bull. Univ. Princeton
(1902), 13.



\bibitem
{Hadamard-1932} J.~Hadamard, Le probl\`{e}me de Cauchy et les
\`{e}quations aux derivees partielles lin\'{e}aires hyperboliques,
Hermann, Paris, 1932.



\bibitem
{Hale} J. K. Hale, "Theory of Functional Differential Equations",
Springer, Berlin- Heidelberg- New York, 1977.


\bibitem
{Hale_book} J. K. Hale and S. M. Verduyn Lunel, "Theory of
Functional Differential Equations", Springer-Verlag, New York, 1993.


\bibitem
{Hartung-Krisztin-Walther-Wu-2006} F.~Hartung, T.~Krisztin,
H.-O.~Walther, J.~Wu, Functional Differential Equations with
State-Dependent Delays: Theory and Applications, {\it in} "Handbook
of Differential Equations: Ordinary Differential Equations,
Volume 3" (A.~Canada, P.~Drabek, A.~Fonda eds.), Elsevier B.V., 2006.

\bibitem
{Hernandez-2006} E.~Hernandez, A.~Prokopczyk, L.~Ladeira,  A note on
partial functional differential equations with state-dependent
delay, Nonlinear Anal. R.W.A. {\bf 7}(4), (2006) 510--519. 

\bibitem
{Lions-Magenes-book}  J.L. Lions and  E. Magenes, 
"Probl\`emes aux Limites Non Homog\'enes et applications". 
 Dunon, Paris, 1968.

\bibitem
{Lions}  J.L.~Lions, "Quelques M\'ethodes de R\'esolution des
Probl\`emes aux Limites Non Lin\'eaires", Dunod, Paris, 1969.

\bibitem%
{Kolmogorov-Fomin} A.N.~Kolmogorov, S.V.~Fomin, "Elements of theory
of functions and functional analysis", Nauka, Moscow, 1968.


\bibitem
{Krisztin-2003} T.~Krisztin, A local unstable manifold for
differential equations with state-dependent delay, Discrete Contin.
Dyn. Syst. {\bf 9}, (2003) 933-1028.



\bibitem%
{Martin-Smith-TAMS-1990} R.H. Martin, Jr., H.L. Smith,  Abstract
functional-differential equations and reaction-diffusion systems,
Trans. Amer. Math. Soc.
{\bf 321} (1990), no. 1, 1--44. 

\bibitem%
{Martin-Smith-JRAM-1991} R.H. Martin, Jr., H.L. Smith,
Reaction-diffusion systems with time delays: monotonicity,
invariance, comparison and convergence, J. reine angew Math. {\bf
413} (1991), 1-35.


\bibitem
{Mallet-Paret} J. Mallet-Paret,  R. D. Nussbaum, P. Paraskevopoulos,
 Periodic solutions for functional-differential equations
 with multiple state-dependent time lags,
  Topol. Methods Nonlinear Anal. {\bf 3}(1), (1994)
  101--162. 



\bibitem
{Rezounenko-Wu-2006} A.V. Rezounenko and J. Wu,  A non-local PDE
model for population dynamics with state-selective delay: local
theory and global attractors, Journal of Computational and Applied
Mathematics, { 190} (1-2), (2006) 99-113.

\bibitem
{Rezounenko_JMAA-2007} A.V. Rezounenko, Partial differential
equations with discrete and distributed state-dependent delays,
Journal of Mathematical Analysis and Applications, {\bf 326}(2),
(2007), 1031-1045. (see also detailed {\it preprint}, March 22,
2005, http://arxiv.org/pdf/math.DS/0503470 ).



\bibitem
{Rezounenko_NA-2009} A.V. Rezounenko, Differential equations with
discrete state-dependent delay: uniqueness and well-posedness in the
space of continuous functions, Nonlinear Analysis Series A: Theory,
Methods and Applications, Volume 70, Issue 11 (2009), 3978-3986.

\bibitem
{Rezounenko-2010} A.V. Rezounenko, Non-linear partial differential
equations with discrete state-dependent delays in a metric space,
Nonlinear Analysis, 73 (2010) 1707-1714 (see also detailed {\it
preprint}: April 15, 2009, http://arxiv.org/pdf/0904.2308v1 ).


\bibitem
{Rezounenko-CR-2011} A.V. Rezounenko, Non-local PDEs with a
state-dependent delay term presented by Stieltjes integral, C. R.
Acad. Sci. Paris, Ser. I 349 (2011), 179-183.

\bibitem
{Ruess-TAMS-2009} W.M. Ruess, Flow invariance for nonlinear partial
differential delay equations,  Trans. Amer. Math. Soc. 361 (2009),
4367-4403.



\bibitem
{So-Wu-Zou} J.~W.-H. So, J.~Wu and X.~Zou, A reaction diffusion
model for a single species with age structure. I. Travelling
wavefronts on unbounded domains, Proc. Royal. Soc. Lond. A { 457},
(2001) 1841-1853.

\bibitem
{So-Yang} J.~W.-H. So and Y.~Yang, Dirichlet problem for the
diffusive    Nicholson's blowflies equation, J. Differential
Equations, {\bf 150}(2), (1998) 317--348. 

\bibitem
{Temam_book} R. Temam, "Infinite Dimensional Dynamical Systems in
Mechanics and Physics", Springer, Berlin-Heidelberg-New York, 1988.

\bibitem
{travis_webb} C.~C. Travis and G.~F. Webb, Existence and stability
for partial functional differential equations,  Transactions of AMS,
{\bf 200},  (1974) 395-418.



\bibitem
{Walther_JDE-2003} H.-O. Walther,  The solution manifold and $C\sp
1$-smoothness for differential equations with state-dependent
delay, J. Differential Equations, {\bf 195}(1),  (2003) 46--65. 


\bibitem
{Walther_JDDE-2007} H.-O. Walther, On a model for soft landing with
state-dependent delay, J. Dynamics and Differential Eqs, {\bf
19}(3), (2007) 593-622.


\bibitem
{Winston-1970} E.~Winston, Uniqueness of the zero solution for
differential equations with state-dependence, J. Differential
Equations, {\bf 7},  (1970) 395--405.

\bibitem
{Wu_book} J.~Wu, "Theory and Applications of Partial Functional
Differential Equations", Springer-Verlag, New York, 1996. 



\bibitem
{yosida} K.~Yosida, "Functional analysis", Springer-Verlag, New
York, 1965.

\end{thebibliography}
\end{document}